
 

\magnification=1200
\input harvmac 
\input amssym.def

\parskip 7pt 
\hfuzz=15pt 
\hoffset -5mm 

\font\male=cmr10

\font\sfont=cmbx10 scaled\magstep2 

\global\newcount\subsecno \global\subsecno=0
\global\newcount\meqno \global\meqno=1
\global\newcount\subsubsecno \global\subsubsecno=0

\def\newsubsec#1{\global\advance\subsecno
by1\message{(\the\secno.\the\subsecno. #1)}
\global\subsubsecno=0 
\noindent{\bf\the\secno.\the\subsecno. #1}\writetoca{\string\quad
{\the\secno.\the\subsecno.} {#1}}}


\def\bu{$\bullet$} 

\def\ssll{\frak{sl}}  
\def\({\left(}
\def\){\right)} \def\cg{{\cal G}} 
\def\cc{{\cal F}} \def\ca{{\cal A}} \def\cf{{\cal F}} 
\def\cb{{\cal B}} \def\bb{{\tilde{\cal B}}} 
 
\def\C{~$\cc$~}  
\def\q{{\tilde q}} \def\tN{{\tilde N}} 
\def\frac#1#2{{#1\over #2}}
 \def\ch{{\cal H}}

\def\vv{\vert n,c\rangle}
\def\vo{~$\vert n,c\rangle$~}
 
\def\eps{\epsilon} 
\def\ww{\langle n,c \vert }
\def\wo{~$\langle n,c \vert $~}
\def\tw{{\tilde w}} \def\tu{{\tilde u}} \def\tL{{\tilde L}} 
\def\tv{{\tilde v}} \def\tp{{\tilde p}} 
\def\hv{{\hat v}} \def\hV{{\hat V}} 
\def\mu{{M}} 
\def\vmr{\vert \mu\rangle}
\def\vmro{~$\vert \mu\rangle$~}

\def\np{\vfil\eject} 
\def\nl{\hfil\break} 
\def\nt{\noindent}

 \def\bbz{Z\!\!\!Z}
\def\bbc{{C\kern-6.5pt I}} \def\bbn{I\!\!N} 
 \def\bbr{I\!\!R} 
  
\def\d{\delta}

\def\l{\lambda}  \def\bU{{\bar U}}
\def\L{\Lambda}

\def\om{\omega} 


\lref\Ji{ M. Jimbo, A $q$-difference analogue of $U(\cg)$ and the
Yang-Baxter equation, Lett. Math. Phys. {\bf 10} (1985) 63-69;\ 
A $q$-difference analogue of $U(gl(N+1))$, Hecke algebras and the
Yang-Baxter equation, Lett. Math. Phys. {\bf 11} (1986) 247-252.} 

\lref\Dr{V.G. Drinfeld, Hopf algebras and the quantum Yang-Baxter
equation, Dokl. Akad. Nauk SSSR {\bf 283} (1985) 1060-1064 (in
Russian); English translation: Soviet. Math. Dokl. {\bf 32}
(1985) 254-258;\ Quantum groups, in: Proceedings ICM 1986, (MSRI,
Berkeley, 1987) pp. 798-820.}

\lref\KRa{P.P. Kulish and N.Yu. Reshetikhin, The quantum linear
problem for the sine-Gordon equation and higher representations,
Zap. Nauch. Semin. LOMI {\bf 101} (1981) 101-110 (in Russian);
English translation: J. Soviet. Math. {\bf 23} (1983) 2435-2441.}
 
\lref\Skc{E.K. Sklyanin, On an algebra generated by quadratic
relations, Uspekhi Mat. Nauk {\bf 40} 
(1985) 214 (in Russian).}

\lref\Sha{ 
N.N. Shapovalov, On a bilinear form on the universal enveloping
algebra of a complex semisimple Lie algebra, Funkts. Anal.
Prilozh. {\bf 6} (4) (1972) 65--70; English translation: Funkt.
Anal. Appl. {\bf 6} (1972) 307--312. }

\lref\Maj{S. Majid, Quantum and braided Lie algebras, { J.
Geom. Phys.} {\bf 13} (1994) 307--356.}  

\lref\DelHuff{G.W. Delius and A. H\"{u}ffmann, On quantum Lie
algebras and quantum root systems, J. Phys. A: Math. Gen. {\bf
29} (1996) 1703-1722.} 

\lref\Del{G.W. Delius, A. H\"uffmann, M. D. Gould and Y.-Z.
Zhang, Quantum Lie algebras associated to $U_q(gl_n)$ and
$U_q(sl_n)$, J. Phys. A: Math. Gen. {\bf
29} (1996) 5611-5617.} 

\lref\BGG{I.N. Bernstein, I.M. Gel'fand and S.I. Gel'fand, Structure
of representations generated by highest weight vectors. Funkts.
Anal. Prilozh. {\bf 5} (1), 1--9 (1971); English translation:
Funkt. Anal. Appl. {\bf 5}, 1--8 (1971)}

\lref\Doa{V.K. Dobrev, Multiplet classification of highest weight
modules over quantum universal enveloping algebras : the
$U_q(sl(3,\bbc))$ example, in: Proceedings of the International Group
Theory Conference (St. Andrews, 1989), eds. C.M.
Campbell et al, Vol. 1, London Math. Soc. Lecture
Note Series 159 (Cambridge University Press, 1991) pp. 87-104.} 

\lref\LS{V. Lyubashenko and A. Sudbery, Quantum Lie Algebras of
Type $A_n\,$, q-alg/9510004, J. Math. Phys., to appear.} 

\lref\Sua{A. Sudbery, The quantum orthogonal mystery, 
in: "Quantum Groups: Formalism and Applications", eds. J.
Lukierski et al (Polish Scientific Publishers PWN, 1995), pp.
303--316.}  

\lref\Sub{A. Sudbery, $SU_q(n)$ gauge theory, Phys. Lett. 
{\bf 375B} (1996) 75-80.} 

\lref\Suc{A. Sudbery, Quantum-group gauge theory, in: "Quantum
Group Symposium at Group21", Proceedings of a Symposium 
at the XXI Intern. Colloquium on Group
Theoretical Methods in Physics, (Goslar, July 1996), eds. H.-D.
Doebner et al (Heron Press, Sofia, 1997) pp. 45-52.} 


\null 

\vskip 1truecm

\centerline{{\sfont Representations of the 
Generalized Lie Algebra ~sl(2)$_{\bf q}$} } 

\vskip 1.5cm

\centerline{{\bf V.K. Dobrev}$^{1}$ 
~and ~{\bf A. Sudbery}$^{2}$}

\vskip 0.75cm

\centerline{$^1$Bulgarian Academy of Sciences} 
\centerline{Institute of Nuclear Research and Nuclear Energy} 
\centerline{72 Tsarigradsko Chaussee, 1784 Sofia, Bulgaria}
\centerline{dobrev@bgearn.acad.bg} 

\vskip 0.5cm 

\centerline{$^2$ University of York}
\centerline{Department of Mathematics} 
\centerline{Heslington, York, England YO1 5DD} 
\centerline{as2@york.ac.uk} 

\vskip 1cm

\centerline{{\bf Abstract}}

\midinsert\narrower\narrower{\male 
We construct finite-dimensional irreducible representations of
two quantum algebras related to the generalized Lie algebra
$\ssll (2)_q\,$ introduced by Lyubashenko and the second named
author. We consider separately the cases of $q$ generic and $q$
at roots of unity. Some of the representations have no classical
analog even for generic $q$. Some of the representations have no
analog to the finite-dimensional representations of the 
quantised enveloping algebra $U_q(sl(2))$, while in those that do
there are different matrix elements.} 
\endinsert

\vskip 10mm
\newsec{Introduction} 

\nt	A number of authors \Maj, \DelHuff, \Del, \LS\ 
have suggested definitions of 
"quantum Lie algebras", the aim being to obtain structures which
bear the same relation to quantised enveloping algebras as Lie
algebras do to their enveloping algebras. It is of interest to
determine the representations of such quantum Lie algebras, in
those cases where a notion of "representation" is defined, and
compare them to the classical representation theory. For generic
values of the deformation parameter $q$ it is to be
expected that the representations will resemble those of the
classical Lie algebras which are deformed into the quantum
versions, since the representation theory of a quantised enveloping
algebra is essentially the same as that of the classical Lie
algebra, but the details of this resemblance will help to
illuminate the nature of a quantum Lie algebra. This relationship
breaks down if $q$ is a root of unity, which is of much interest
in physics, and it is therefore particularly significant to
determine the representations of a quantum Lie algebra in this
case. 

In this paper we start on such a study by constructing 
finite-dimensional representations of the simplest example of the
generalized Lie algebras introduced in \LS. A representation of
this algebra, in the sense defined in \LS, is nothing but a
representation of an associative algebra, the enveloping algebra
of the quantum Lie algebra. This is obtained from a larger
algebra with a central element by imposing a relation giving the
central element as a function of Casimir-like elements. We
investigate the representations also of this larger algebra,
which is possibly more natural in the context of generalized Lie
algebras, and find that it has additional one-dimensional
representations.

The paper is organised as follows. In Section 2 we introduce
explicitly the two quantum algebras which we consider. In
Sections 3 and 4 we construct finite-dimensional representations
of these algebras for generic values of $q$. In Sections 5 and 6
we consider the case when $q$ is at roots of unity. Section 7
contains a Summary of our results. 

\vskip 10mm 

\newsec{The quantum Lie algebra $\ssll (2)_q$} 

\nt
The generalized Lie algebra ~$\ssll (2)_q$ ~was introduced in
\LS, cf. also \Sua, \Sub, \Suc. Its enveloping algebra ~$\ca
\equiv U(\ssll (2)_q)$ ~ is defined by Eq. (3.5) of \LS. For the
purposes of developing the representation theory it is enough to
work with the algebras $\cb$, $\cc$, cf. \LS. The algebra ~$\cb$~
is generated by four generators:~ $X_0,X_\pm,C$ ~with relations: 
\eqna\XC $$\eqalignno{ 
& q^2X_0X_+-X_+X_0 ~=~qCX_+ &\XC a\cr 
& q^{-2}X_0X_--X_-X_0 ~=~- q^{-1}CX_- &\XC b\cr 
& X_+X_--X_-X_+ ~=~ 
(q+q^{-1})\ (C-\l X_0) X_0 &\XC c\cr 
& CX_m ~=~ X_mC ~,\qquad \qquad m=0,\pm 1 &\XC d\cr 
& \qquad \l \equiv q -q^{-1} }$$   

The algebra $\cb$ is related to the locally finite part $\cf$ of
the simply-connected quantised enveloping algebra $\bU_q(sl(2))$.
The algebra $\cf$ was obtained in \LS\ from $\cb$ by putting 
$C$ equal to a function of the second order Casimir: 
\eqn\cas{ C_2 ~~=~~ (q+q^{-1})\ X_0^2 ~+~ qX_-X_+ ~+~
q^{-1}X_+X_- } 
namely, 
\eqn\css{ C^2 ~~=~~ 1 ~+~ {\l^2 \over q+q^{-1}}\ C_2 } 
For shortness we shall call ~$\cc$~ the restricted algebra.  
The enveloping algebra $\ca$, on the other hand, is obtained by
putting $C=1$ \LS. 

We shall need a triangular decomposition of $\cb$~: 
\eqn\tri{\cb ~~=~~ \cb_+ \otimes \cb_0 \otimes \cb_- }
where $\cb_\pm$ is generated by $X_\pm$, while $\cb_0$ is generated
by $X_0,C$. We shall call the abelian Lie algebra $\ch$
generated by $X_0,C$ the Cartan subalgebra of $\cb$. Note that
$\cb_0$ is the enveloping algebra of $\ch$.  The same decomposition 
is used for the algebra $\cc$ with the relation \css\ enforced. 

Further we shall analyse the algebras ~$\cb$~ and ~$\cc$~ separately. 
 
\vskip 10mm 

\newsec{Highest weight representations} 

\nt 
Highest weight modules of $\cb$ are standardly determined
by a highest weight vector $v_0$ which is annihilated by the
raising generator $X_+$ and on which the Cartan generators act
by the corresponding value of the highest weight $\L\in\ch^*$~:
\eqn\hwm{ \eqalign{ 
X_+\ v_0 ~~=&~~ 0 \cr
H\ v_0 ~~=&~~ \L(H)\ v_0 ~, \qquad H\in\ch \cr
&\mu \equiv \L(X_0) ~, \quad c \equiv \L(C) \cr }} 

In particular, we shall be interested in Verma modules over \C.
As in the classical case a Verma module ~$V^\L$~ is a highest 
weight module (HWM) of weight $\L$ induced from one-dimensional
representation of a Borel subalgebra $\bb$, e.g., $\bb = 
\cb_+ \otimes \cb_0\,$, on the highest weight vector, e.g.,
$v_0$. As vector spaces we have: 
\eqn\verm{ 
V^\L ~\cong ~ \cb \otimes_\bb v_0 ~=~ \cb_- \otimes
v_0 ~=~ l.s. \{ v_k \equiv X_-^k \otimes v_0 ~\vert~ k\in\bbz_+
\}} 
where we have identified ~$1_\cb \otimes v_0$~ with ~$v_0\,$. 
 
The action of the generators of $\cb$ on the basis of $V^\L$ is
given as follows: 
\eqna\act
$$\eqalignno{
X_+\ v_k ~~=&~~ q^{2k-2}\ \left(c-\l \mu\right)\ 
\left( [2k]_q \mu -q [k]_q [k-1]_q c \right) \ v_{k-1} &\act a\cr 
X_-\ v_k ~~=&~~ v_{k+1} &\act b\cr 
X_0\ v_k ~~=&~~ ( q^{2k} \mu - q^k [k]_q c )\ 
v_{k} &\act c\cr 
C\ v_k ~~=&~~ c\ v_k &\act d\cr 
&[k]_q ~\equiv~ (q^k-q^{-k})/\l 
}$$ 
To obtain \act{a,c} we have used the following calculations which
follow from \XC{}~: 
\eqna\clc 
$$\eqalignno{ 
X_0 X_-^k ~~=&~~ X_-^k \left( q^{2k} X_0 - q^k [k]_q C
\right)  &\clc a\cr 
[X_+, X_-^k] ~~=&~~ X_-^{k-1}\ q^{2k-2}\ \left(C-\l X_0\right)\ 
\left( [2k]_q X_0 -q [k]_q [k-1]_q C \right) &\clc b\cr 
}$$ 

As in the classical case the analysis of reducibility of Verma
modules is an important tool in the representation theory. This
analysis starts (cf. \BGG) with the search for singular vectors. 
A singular vector ~$v_s$~ of a Verma
module ~$V^\L$~ is defined as follows: ~$v_s \in V^\L$~, ~$ v_s
\notin \bbc v_0$ and it satisfies the following properties (cf.,
e.g., \BGG)~: 
\eqna\sing 
$$\eqalignno{ 
X_+\ v_s ~~=&~~ 0 &\sing a\cr
H \ v_s ~~=&~~ \L'(H)\ v_s ~, \qquad H\in\ch ~, ~~~ 
\L'\in\ch^* &\sing b\cr}$$ 
{}First we note that since $C$ is central its value is the same
as on $v_0$~: ~$c' \equiv \L'(C) = c$. Further, we proceed to
find the possible singular vectors using that they are
eigenvectors of $X_0\,$. But the eigenvectors of ~$X_0$~ are
~$X_-^n\otimes v_0\,$, all with different eigenvalues. 
Thus, a singular vector will be 
given by the classical expression (omitting
the overall normalization)~: ~$v_s ~=~ 
X_-^n \otimes v_0$~ for some fixed $n\in\bbn$ , and we have:
\eqn\vgh{ X_0\ v_s ~=~ \mu'\ v_s ~, \qquad 
\mu' ~\equiv ~ \L'(X_0) ~=~ q^{2n} \mu - q^n [n]_q c } 
Finally, we have to impose \sing{a} for which we calculate (using
\clc{b})~:
\eqn\vgs{ X_+\ v_s ~~=~~ X_-^{n-1}\ q^{2n-2}\ \left(c-\l \mu\right)\ 
\left( [2n]_q \mu -q [n]_q [n-1]_q c \right) \otimes v_0 } 
For the further analysis we suppose that the deformation
parameter $q$ is not a nontrivial root of unity. 
Then there are two possibilities for \vgs\ to be zero, and 
thus, we have two possibilities to fulfil \sing{a}~:
\eqna\vgp
$$\eqalignno{ \mu ~=&~ q [n]_q [n-1]_q c / [2n]_q &\vgp a\cr 
c ~=&~ \l\mu &\vgp b\cr }$$ 

We shall analyse the two possibilities in \vgp{} separately since
they have very different implications; moreover, they are
incompatible unless $c=M=0$ when they coincide and which we shall
treat as partial case of \vgp{b}. 

\newsubsec{} The first possibility \vgp{a} (with $c\neq 0$)
corresponds to the classical 
relation between $n$ and the highest weight $\L$ (obtained for
$q,c\to 1$)~: ~$\mu ~=~ (n-1)/2$. Thus, if we fix ~$n\in\bbn$~
then ~$v_s ~=~ X_-^n \otimes v_0$~ is a singular vector when
~$\mu$~ has the value \vgp{a}. The shifted weight ~$\L'$~
corresponds to another Verma module ~$V^{\L'}$~ which is the 
naximal invariant submodule of ~$V^\L$. The corresponding
eigenvalue of $X_0$ is (cf. \vgh): 
\eqn\mmm{ \mu' ~=~ - q [n]_q [n+1]_q c / [2n]_q } 
Note that the Verma module $V^{\L'}$ does not have a singular
vector. Indeed, there is no ~$n'\in\bbn$~ such that \vgp{a} holds
for the pair ~$(\mu',n')$~ replacing $(\mu,n)$. Also \vgp{b} can
not hold for $\mu'$ since $c=\l\mu'$ will contradict \mmm. 

The factor-module ~$L_{n,c} \cong V^\L/V^{\L'}$~ is irreducible and 
finite-dimensional of dimension ~$n$. It has a highest weight
vector \vo such that: 
\eqn\hmm{ \eqalign{ 
X_+\ \vv ~~=&~~ 0 \cr
H\ \vv ~~=&~~ \L(H)\ \vv ~, \qquad H\in\ch \cr
X_-^n\ \vv ~~=&~~ 0 \cr }} 
Let us denote by ~$w_k \equiv X_-^{k}\vv$, ~$k=0,1,\dots,n-1$,
the states of $L_{n,c}\,$. The transformation rules for ~$w_k$~ are:
\eqna\fct 
$$\eqalignno{
X_+\ w_k ~~=&~~ q^{2k-n}\ [k]_q [n-k]_q \left( {c\ [2]_q\ [n]_q 
\over [2n]_q }\right)^2 \ w_{k-1} &\fct a\cr 
X_-\ w_k ~~=&~~ w_{k+1} ~, \qquad k<n-1 &\fct b\cr 
X_-\ w_{n-1} ~~=&~~ 0 &\fct {b'}\cr 
X_0\ w_k ~~=&~~ {c\ q^k\ [n]_q \over [2n]_q }\  
 \left( [n-k]_q - q^{1-n} [k+1]_q \right) \ 
w_{k} &\fct c\cr 
C\ w_k ~~=&~~ c\ w_k &\fct d\cr 
}$$ 
Thus, the vector ~$w_{n-1}$~ is the lowest weight vector of
$L_{n,c}\,$. 

Next we introduce a bilinear form in $L_{n,c}$ by the formula:
\eqn\sca{ (w_j,w_k) ~\equiv~ \ww X^j_+ X^k_- \vv }
where \wo is such that ~$\ww\ \vv ~=~ 1$~ and:
\eqn\hml{ \eqalign{ \ww\ X_- ~~=&~~ 0 \cr 
\ww\ H ~~=&~~ \L(H)\ \ww ~, \qquad H\in\ch \cr
\ww\ X_+^n ~~=&~~ 0 \cr }} 
Then we obtain~:
\eqn\scb{ \eqalign{ (w_j,w_k) ~=&~ \d_{jk}\ 
q^{k(k+1-n)}\ { [k]_q!\ [n-1]_q! 
\over [n-1-k]_q! }\ \left( {c\ [2]_q\ [n]_q \over [2n]_q
}\right)^{2k} \cr 
&[k]_q! ~\equiv ~[k]_q[k-1]_q\dots [1]_q ~, 
\quad [0]_q! ~\equiv ~ 1 }}
Clearly, \scb\ is real-valued for real ~$q,c$. Thus, for $q,c\in\bbr$
we can turn \sca\ into a scalar product and define the norm
of the basis vectors: 
\eqn\nor{ \vert w_k\vert ~\equiv ~ \sqrt{(w_k,w_k)} ~=~ 
q^{k(k+1-n)/2}\ \sqrt{{ [k]_q!\ [n-1]_q! 
\over [n-1-k]_q! }}\ 
\left( {c\ [2]_q\ [n]_q \over [2n]_q }\right)^{k} } 
where have chosen the root that is positive for positive $c,q$. 
We can also introduce orthonormal basis:
\eqn\ort{ u_k ~\equiv {1\over \vert w_k\vert}\ w_k } 
Then we have: 
\eqn\scc{ (u_j,u_k) ~=~ \d_{jk} }

The transformation rules for the basis vectors ~$u_k$~ are:
\eqna\uct 
$$\eqalignno{
X_+\ u_k ~~=&~~ q^{k-n/2}\ \sqrt{[k]_q [n-k]_q}\  
 {c\ [2]_q\ [n]_q \over [2n]_q } \ u_{k-1} &\uct a\cr  
X_-\ u_k ~~=&~~ q^{k+1-n/2}\ \sqrt{[n-1-k]_q \ [k+1]_q}\  
 {c\ [2]_q\ [n]_q \over [2n]_q }\ u_{k+1} &\uct b\cr 
X_0\ u_k ~~=&~~ {c\ q^k\ [n]_q \over [2n]_q }\  
 \left( [n-k]_q - q^{1-n} [k+1]_q \right)\ 
u_{k} &\uct c\cr 
C\ u_k ~~=&~~ c\ u_k &\uct d\cr 
}$$ 

The above scalar product is invariant under the real form
~$\cb_r$~ of ~$\cb$~ defined by the antilinear antiinvolution:
\eqn\cnj{ \om (X^\pm) ~=~ X^\mp ~, \qquad \om (X_0) ~=~ X_0 
~, \qquad \om (C) ~=~ C } 
Indeed, the algebraic relations \XC{} are preserved by ~$\om$~
for real ~$q$. The ~$\cb_r$~ invariance of the scalar product
means that:
\eqn\scd{ (w_j,X w_k) ~=~ (\om(X) w_j,w_k) ~, \qquad X\in\cb, } 
which is automatically satisfied with the definition \sca. 
(Note that \scd\ defines ~$(,)$~ as the Shapovalov bilinear form
\Sha.) 

Thus, for every ~$n\in\bbn$~ we have constructed ~$n$-dimensional
irreducible representations (irreps) of ~$\cb$~ parametrized by
~$c\in\bbc$, $c\neq 0$, with basis ~$w_k$~ or ~$u_k\,$, 
($k=0,...,n-1$). For ~$q,c\in\bbr$~ these are irreps of the real
form ~$\cb_r\,$, which are unitary when ~$q,c>0$.

\newsubsec{} 
The second possibility \vgp{b} has no classical analogue. It
tells us that if ~$c$~ and ~$\mu$~ are related as in \vgp{b} then
each vector of the basis of $V^\L$ is a singular vector.
Moreover, all of them have the same weight since ~$\mu' ~=~ \mu$,
cf. \vgh. This is clear also from the transformation rules \act{}
when ~$c=\l\mu$~:
\eqna\mct $$\eqalignno{
X_+\ v_k ~~=&~~ 0 &\mct a\cr 
X_-\ v_k ~~=&~~ v_{k+1} &\mct b\cr 
X_0\ v_k ~~=&~~ \mu\ v_{k} &\mct c\cr 
C\ v_k ~~=&~~ \l\mu\ v_k &\mct d\cr }$$ 
Clearly, we have an infinite sequence of embedded reducible Verma 
modules\nl $V_n ~=~ l.s. \{ v_k ~\vert~ k\in\bbz_+ ~, \quad k\geq
n\}$ for ~$n\in\bbz_+\,$~ as follows: ~$V_n ~\supset~ V_{n+1}\,$,
the latter being the maximal invariant submodule of the former. 
Note that ~$V_n$~ is isomorphic to a submodule of all ~$V_m$~
with ~$n>m$. Furthermore, because of the coincidence of the
weights these modules are also all isomorphic to each other:
~$V_n \cong V_m$~ for all ~$m,n$. It is also clear that for every
$\mu$ there is only one irreducible module, namely the
one-dimensional ~$L_\mu ~\cong V_n/V_{n+1}\,$, for any $n$.
Denoting by \vmro the only state in ~$L_\mu$~ we have for the
action on it: 
\eqna\nct
$$\eqalignno{
X_+\ \vmr ~~=&~~ 0 &\nct a\cr 
X_-\ \vmr ~~=&~~ 0 &\nct b\cr 
X_0\ \vmr ~~=&~~ \mu\ \vmr &\nct c\cr 
C\ \vmr ~~=&~~ \l\mu\ \vmr &\nct d\cr 
}$$ 

Note that the above one-dimensional irrep is different from the
one-dimenional ~$L_{1,c}$~ from the previous subsection.
Indeed, though the action of $X_\pm$ is the same, the 
ratio of eigenvalues of $C$ to $X_0$ here is 
$\l$, while there it is $-[2]_q/q$~. 
 
\vskip 10mm 
 
\newsec{Highest weight representations of the restricted algebra} 

\nt 
The highest weight representations of the restricted algebra
~$\cc$~ are obtained from those of ~$\cb$~ imposing the relation
\css. In particular, there is the following relation
between the values of the Cartan generators:
\eqn\crel{ c^2 ~~=~~ 1 + \l^2\ \( \frac{\mu^2}{q^2}
+ c\frac{\mu}{q} \) }
This relation has to be imposed on all formulae of the previous
Section. There are no essential consequences of this for the
generic Verma modules. For the reducible Verma modules there are
more interesting consequences. First we notice that the 
reducibility condition \vgp{b} is incompatible with \crel, and
thus there would be no special one-dimensional irreps like
$L_\mu\,$, cf. \nct. So it remains to consider the combination of the
reducibility condition \vgp{a} with \crel\ from which we obtain that: 
\eqn\vgr{ c ~=~ {\eps\ [2n]_q 
\over [2]_q\ [n]_q }  ~, \qquad \mu ~=~ {q [n]_q [n-1]_q\ c
\over [2n]_q }  ~=~ {\eps\ q [n-1]_q \over [2]_q} 
~, \qquad \eps = \pm 1 } 
In this case the analogue of \mmm\ is: 
\eqn\mmz{ \mu' ~=~ - \eps\ q [n+1]_q  / [2]_q } 
 Let us denote the finite-dimensional
representations of $\cc$ by ~$\tilde L_{n,\eps}$~ and the basis by
~$\tw_k\,$, $k=0,...,n-1$. The transformation rules are:
\eqna\wct 
$$\eqalignno{
X_+\ \tw_k ~~=&~~ q^{2k-n}\ [k]_q [n-k]_q\  \tw_{k-1} &\wct a\cr 
X_-\ \tw_k ~~=&~~ \tw_{k+1} ~, \qquad k<n-1 &\wct b\cr 
X_-\ \tw_{n-1} ~~=&~~ 0 &\wct {b'}\cr 
X_0\ \tw_k ~~=&~~ {\eps\ q^k\over [2]_q }\ 
 \left(  [n-k]_q - q^{1-n} [k+1]_q \right)\ \tw_{k} &\wct c\cr 
C\ \tw_k ~~=&~~ {\eps\ [2n]_q \over [2]_q [n]_q }\ \tw_k &\wct d\cr 
}$$ 
Further, the analogues of \scb\ and \nor\ are:
\eqn\scbb{ \eqalign{ (\tw_j,\tw_k) ~=~ \d_{jk}\ 
q^{k(k+1-n)}\ { [k]_q!\ [n-1]_q! 
\over [n-1-k]_q! } }}
\eqn\norr{ \vert \tw_k\vert ~\equiv ~ \sqrt{(\tw_k,\tw_k)} ~=~ 
q^{k(k+1-n)/2}\ [k]_q!\ \sqrt{{[k]_q!\ [n-1]_q! 
\over [n-1-k]_q! }} } 
We can also introduce orthonormal basis:
\eqn\ortt{ \tu_k ~\equiv {1\over \vert \tw_k\vert}\ \tw_k ~, 
\qquad (\tu_j,\tu_k) ~=~ \d_{jk} } 
for which the transformation rules are: 
\eqna\uctt 
$$\eqalignno{
X_+\ \tu_k ~~=&~~ q^{k-n/2}\ \sqrt{[k]_q [n-k]_q} 
 \ \tu_{k-1} &\uctt
a\cr  
X_-\ \tu_k ~~=&~~ q^{k+1-n/2}\ \sqrt{[n-1-k]_q \ [k+1]_q}\ 
\tu_{k+1} &\uctt b\cr 
X_0\ \tu_k ~~=&~~ {\eps\ q^k\over [2]_q }\ 
 \left( [n-k]_q - q^{1-n} [k+1]_q \right)\ \tu_{k} &\uctt c\cr 
C\ \tu_k ~~=&~~ {\eps\ [2n]_q \over [2]_q [n]_q }\ \tu_k &\uctt d\cr 
}$$ 

Thus, for every ~$n\in\bbn$~ we have constructed ~$n$-dimensional
irreducible representations of ~$\cc$~ parametrized by ~$\eps
=\pm 1$, with bases ~$\tw_k\,$ or $\tu_k\,$, ($k=0,...,n-1$). 

\newsec{Highest weight representations at roots of unity} 

\nt 
Here we consider representations of the algebra $\cb$ in the case
when the deformation parameter is at 
roots of unity. More precisely, first we consider the cases when
~$q^2$~ is a primitive $N$-th root of unity: ~$q ~=~ e^{\pi
i/N}$, $N\in\bbn+1$. Then we have: 
\eqn\rrr{ [x]_q ~=~ {\sin \pi x /N \over \sin \pi /N } } 
In such cases there are additional reducibility conditions
coming from \vgs\ besides \vgp{a,b}. For this we rewrite \vgp{a} in
a more general fashion:
$$\eqalignno{ \mu [2n]_q ~=&~ q [n]_q [n-1]_q c  &\vgp {a'}}$$ 
Then we note that from \rrr\ follows
that ~$[N]_q ~=~[2N]_q ~=~ 0$, so \vgp{a'} is satisfied for
~$n\to N$. Thus, ~$v_s^N ~=~ X_-^N
\otimes v_0$ is a singular vector independently of the highest
weight $\L$. Similarly to the analysis done in \Doa\ for the
quantised enveloping algebra\foot{We recall that
though the quantised enveloping algebras ~$U_q(\cg))$~ were
introduced for arbitrary simple Lie algebras $\cg$ in \Dr, \Ji, 
the example of $U_q(sl(2))$ was introduced in 
\KRa\ as an algebra and in \Skc\ as a Hopf algebra.} ~$U_q(sl(2))$~ 
all ~$v_s^{pN} ~=~ X_-^{pN}\otimes v_0$~ for $p\in\bbn$ are
singular vectors. The Verma modules they realize we denote by
~$\tilde V_p\,$, $p\in\bbz_+\,$, $\tilde V_0 \equiv V^\L$. These
are embedded reducible Verma modules ~$\tilde V_p \supset \tilde
V_{p+1}$~ with the same highest weight $\L$. Indeed, for any
~$\tilde V_p$~ using \vgh\ with $n\to pN$ we have: ~$\mu' ~=~
q^{2pN} \mu - q^{pN} [pN]_q c ~=~ \mu$.

The further analysis depends on whether there are additional singular 
vectors besides those just displayed. There are four cases. 

\newsubsec{} We start with the case when ~$\mu,c$~ do not satisfy
either of \vgp{a,b}. We also suppose that ~$c\neq 0$ when ~$N$~ is ~{\it
even}. Then there are no additional singular vectors and there 
is only one irreducible $N$-dimensional HWM ~$L_{\L,N} 
\cong \tilde V_p/\tilde V_{p+1}$~ (for any $p$), parametrized by
all pairs ~$\mu,c$~ not satisfying \vgp{a,b}. The action of the
generators of $\cb$ on the basis of ~$L_{\L,N}\,$, which we denote by
~$\tv_k\,$, ($k=0,...,N-1$), is given as follows: 
\eqna\vct
$$\eqalignno{
X_+\ \tv_k ~~=&~~ q^{2k-2}\ \left(c-\l \mu\right)\ 
\left( [2k]_q \mu -q [k]_q [k-1]_q c \right) \ \tv_{k-1} &\vct a\cr 
X_-\ \tv_k ~~=&~~ \tv_{k+1} ~, \qquad k<N-1
&\vct b\cr 
X_-\ \tv_{N-1} ~~=&~~ 0 &\vct {b'}\cr 
X_0\ \tv_k ~~=&~~ ( q^{2k} \mu - q^k [k]_q c )\ 
\tv_{k} &\vct c\cr 
C\ \tv_k ~~=&~~ c\ \tv_k &\vct d\cr }$$ 
However, unlike the D-J case, these finite-dimensional
representations are not unitarizable, which is easily seen if one
considers the analogue of the bilinear form \sca. 

\newsubsec{} 
Next we consider the case when ~$\mu,c$~ satisfy \vgp{a} 
for some ~$n\in\bbn$, $n<N$. We also suppose that ~$c\neq 0$ (for
any $N$). First we note that ~$n<N$~ is not a
restriction, since then \vgp{a} holds also 
for all ~~$n+pN$, $p\in\bbz$. Indeed, we have:
\eqn\sss{\eqalign{ 
&~ q [n+pN]_q [n+pN-1]_q c / [2(n+pN)]_q ~=~ 
q [n]_q [n-1]_q \cos^2 (\pi p)\ c / [2n]_q \cos (2\pi p) ~=\cr 
&=~ q [n]_q [n-1]_q c / [2n]_q ~=~ \mu }} 
Thus, we have another infinite series of singular vectors ~
~$v'^{pN}_s ~=~ X_-^{n+pN}\otimes v_0$~ for $p\in\bbz_+\,$.
They realize reducible Verma modules which we denote by
~$\tilde V'_p\,$, $p\in\bbz_+\,$; ($\tilde V'_0$~ is the analogue
of ~$V^{\L'}$~ introduced in the non-root-of-unity case, but here
it is reducible). 
They all have the same highest weight ~$\L'$ determined by $\mu',c$
with $\mu'$ given by \vgh. Indeed, substituting $n$ with $n+pN$
does not change the value of $\mu'$~: 
\eqn\ssa{\eqalign{ 
&q^{2(n+pN)} \mu - q^{n+pN} [n+pN]_q c ~=~ 
q^{2n} \mu - q^{n+pN}\ e^{\pi i p}\ [n]_q\ \cos (\pi p)\ c ~=\cr 
&=~ q^{2n} \mu - q^{n} [n]_q c ~=~ \mu' }} 
Of course, after substituting $\mu$ with its value from \vgp{a}
we obtain the expression for $\mu'$ in \mmm. 
We have the following infinite embedding chain:
\eqn\emb{ V^\L\equiv \tilde V_0 \supset \tilde V'_0 
\supset \tilde V_1 \supset \tilde V'_1 \supset ... } 
where all embeddings are non-composite: ~
the embeddings ~$\tilde V_p \supset \tilde V'_{p}$~ are
realized by singular vectors: ~$X_-^{n}\otimes v_p\,$,~ 
$v_p$ being the highest weight vector of $\tilde V_p\,$, while 
the embeddings ~$\tilde V'_p \supset \tilde V_{p+1}$~ are
realized by singular vectors: ~$X_-^{N-n}\otimes v'_p\,$,~ 
$v'_p$ being the highest weight vector of $\tilde V'_p\,$. 

Now, factorizing each reducible Verma module by its maximal
invariant submodule we obtain that for each ~$n\in\bbn$, $n<N$~
there are two finite dimensional irreps parametrized by
$c\in\bbc$, $c\neq 0$~: ~$L_{n,N} 
\cong \tilde V_p/\tilde V'_{p}$~ (for any $p$) which
is ~$n$-dimensional, and ~$L'_{n,N} 
\cong \tilde V'_p/\tilde V_{p+1}$~ (for any $p$) which
is ~$(N-n)$-dimensional. However, it turns out that the irreps
from one series are isomorphic to those of the other: ~$L'_{n,N} 
\cong L_{N-n,N}\,$. Indeed, note that the value of ~$\mu'$~ for 
the Verma modules ~$\tilde V'_p$~ given by \mmm\ should be
obtained (for consistency) also from the formula for $\mu$ with
$n$ substituted by $N-n$ (since this is the reducibility
condition w.r.t. the non-composite singular vector 
~$X_-^{N-n}\otimes v'_p\,$) and indeed this is the case:
$$\eqalign{ &q [N-n]_q [N-n-1]_q c / [2(N-n)]_q ~=~ - 
q [n-N]_q [n+1-N]_q c / [2(n-N)]_q ~=\cr &= ~- 
q [n]_q [n+1]_q c\ \cos^2 (\pi N) / [2n]_q\ \cos (2\pi N) ~=~ 
- q [n]_q [n+1]_q c / [2n]_q ~=~ \mu' } $$ 
Furthermore, the transformation rules for 
~$L_{n,N}$~ are the same as for ~$L_{n,c}\,$, cf. \fct{},  
while the transformation rules for $L'_{n,N}$ are obtained 
from \fct{} by substituting $n\to N-n$. 

Thus, we are left with one series of finite-dimensional irreps
~$L_{n,N}\,$. 

\newsubsec{} 
Next, we consider the case when ~$\mu,c$~ satisfy \vgp{b} 
for arbitrary ~$c$. Actually, nothing is changed from the
non-root-of-unity case since the relevant formulae 
\mct{} and \nct{} are not changed. 

\newsubsec{} 
Finally, we consider the case when $N$ is even and $c=0$. Let
~$\tN = N/2 \in\bbn$. In these cases there are 
additional reducibility conditions coming from \vgp{a'}. 
Indeed, from \rrr\ follows that ~$[2\tN]_\q ~=~ 0$~ and
~$[\tN]_\q ~\neq 0$. But if ~$c=0$~ then \vgp{a'} is again
satisfied. Thus, 
the vector ~$\hv_s^{\tN} ~=~ X_-^{\tN} \otimes v_0$ is a singular vector
independently of the value of $\mu$. Similarly to the analysis of
the first subsection also all 
~$\hv_s^{p\tN} ~=~ X_-^{p\tN}\otimes v_0$~ for $p\in\bbn$ are
singular vectors. Note that for $p$ even these are the 
singular vectors that we already have: ~$\hv_s^{p\tN} = v_s^{\tp N}$, 
$\tp =p/2$. The Verma modules they realize we denote by
~$\hat V_p\,$, $p\in\bbz_+\,$, $\hat V_0 \equiv V^\L$. These are embedded
reducible Verma modules ~$\hat V_p \supset \hat V_{p+1}$~ with the same
value of $\mu$ up to sign. Indeed, for any 
~$\hat V_p$~ using \vgh\ with $n\to p\tN$ we have: ~$\mu' ~=~
q^{2p\tN} \mu - q^{p\tN} [p\tN]_q c ~=~ (-1)^p\mu$.
Certainly, for even $p$ these are Verma modules from the first 
subsection: ~$\hV_{p} = V_{\tp}$. 

As above the further analysis depends on whether ~$\mu,c$~
satisfy some of \vgp{a,b}. However, since $c=0$ then the only
additional possibility is that also $\mu=0$, which is a partial
case of \vgp{b}, which was considered in the previous subsection.
Thus, further, we suppose that ~$\mu,c$~ do not 
satisfy either of \vgp{a,b} and that $\mu\neq 0$. 

Then there are no additional singular vectors besides $\hv_s^{p\tN}$. 
Then for each ~$\mu\neq 0$~ 
there is only one irreducible HWM ~$ L_{\mu,\tN} 
\cong \hat V_p/\hat V_{p+1}$~ (for any $p$) which
is ~$\tN$-dimensional. The action of the generators of $\cb$ on the
basis of ~$ L_{\mu,\tN}\,$, which we denote by ~$\hv_k\,$,
($k=0,...,\tN-1$), is given as follows: 
\eqna\vctt
$$\eqalignno{
X_+\ \hv_k ~~=&~~ - q^{2k-2} \l 
[2k]_q \mu^2 \ \hv_{k-1} &\vctt a\cr 
X_-\ \hv_k ~~=&~~ \hv_{k+1} ~, \qquad k<\tN-1
&\vctt b\cr 
X_-\ \hv_{\tN-1} ~~=&~~ 0 &\vctt {b'}\cr 
X_0\ \hv_k ~~=&~~ q^{2k} \mu \ 
\hv_{k} &\vctt c\cr 
C\ \hv_k ~~=&~~ 0 &\vctt d\cr }$$ 
Note that if ~$\tN$~ is odd it seems that formulae \vctt{} may be
obtained from \vct{} for ~$N$ odd and $c=0$ by the substitution
$N\to \tN$. However, this is not the same irrep since with the
same replacement the parameter $q$ there becomes $e^{\pi i/N} \to
e^{\pi i/\tN}$ while the parameter $q$ here is ~$e^{\pi i/2\tN}$.

\newsec{Highest weight representations at roots of unity of the
restricted algebra} 

\nt 
Here we consider representations of the restricted algebra $\cc$
in the case when the deformation parameter is at roots of unity.
We start with the case: ~$q ~=~ e^{\pi i/N}$, $N\in\bbn+1$, and
so \vgp{a'} holds. The analysis is as for the algebra $\cb$ but
imposing the relation \crel, i.e., combining the considerations
of the previous two Sections. 

\newsubsec{} 
We start with the case when ~$\mu,c$~ do not satisfy 
\vgp{a}, i.e., \vgr\ does not hold. We also suppose that ~$c\neq 0$ when
~$N$~ is ~{\it even}. Then there is only one irreducible
$N$-dimensional HWM parametrized by $\mu,c$ related by \crel, 
which irrep we denote by ~$\tL_{\L,N}\,$. For the transformation
rules we can use formulae \vct{} with \crel\ imposed.  

\newsubsec{} 
Next we consider the case when ~$\mu,c$~ satisfy 
\vgp{a} and $c\neq 0$. Here we should be nore careful so we
replace ~$n$~ by ~$n+pN$~ with $n<N$. Combining the 
reducibility condition \vgp{a} with \crel\ we first obtain that: 
\eqn\vgrz{ c^2 ~=~ { [2(n+pN)]_q^2 \over [2]_q^2\ [n+pN]_q^2 } 
~=~ { [2n]_q^2 \over [2]_q^2\ [n]_q^2 } } 
Then we recover \vgr\ and \mmz\ for $n<N$ which means that we have 
the same situation as for the unrestricted algebra at roots of 
unity. Thus, for each ~$n\in\bbn$, $n<N$ and $\eps=\pm1$
there is a finite dimensional irrep: ~$\tL_{n,\eps,N}$~ 
which is ~$n$-dimensional. The transformation rules for 
~$\tL_{n,\eps,N}$~ are the same as in the non-root-of-unity case,
cf. \wct{}. 

\newsubsec{} 
Finally, we consider the case when $N$ is even and $c=0$. Let
~$\tN = N/2 \in\bbn$. 
As for the unrestricted algebra there are additional 
reducibility conditions, i.e., again the vector 
~$v_s^{\tN} ~=~ X_-^{\tN} \otimes v_0$ is a singular vector. 
However, because of \crel\ the value of $\mu^2$ is fixed:
\eqn\rru{ \mu^2 ~=~ -\q^2/\l^2 ~, \qquad \mu ~=~ \eps\ i \q/\l 
~, \qquad \eps =\pm1 } 
Otherwise, the analysis goes through and 
there is only one irreducible ~$\tN$-dimensional HWM 
~$\tL_{\eps,\tN}$~ parametrized by $\eps$. The action of the
generators of $\cb$ on the basis of ~$\tL_{\eps,\tN}\,$, which we denote by
~$\hv'_k\,$, ($k=0,...,\tN-1$), is given as follows: 
\eqna\vctt
$$\eqalignno{
X_+\ \hv'_k ~~=&~~ {\q^{2k} 
[2k]_\q \over \l} \ \hv'_{k-1} &\vctt a\cr 
X_-\ \hv'_k ~~=&~~ \hv'_{k+1} ~, \qquad k<\tN-1
&\vctt b\cr 
X_-\ \hv'_{\tN-1} ~~=&~~ 0 &\vctt {b'}\cr 
X_0\ \hv'_k ~~=&~~ {\eps\ i \q^{2k+1} \over\l} \ 
\hv'_{k} &\vctt c\cr 
C\ \hv'_k ~~=&~~ 0 &\vctt d\cr }$$ 
The crucial feature of these two irreps is that they do not have
a classical limit for $\q\to 1$ (obtained for $N\to\infty$).

\newsec{Summary} 

Below by ~$q$~ {\it generic}~ we shall understand that $q$ is a
nonzero complex number which is not a nontrivial root of unity. 
We have constructed the following finite-dimensional irreps of
the algebras ~$\cb$~ and ~$\cc$. 

\newsubsec{} 
For the algebra ~$\cb$~: 

\item{\bu} $L_{n,c}\,$, ~$n\in\bbn$, ~$c\in\bbc$, $c\neq 0$,
~$q$~ generic, ~$\dim L_{n,c} ~=~ n$, cf. \fct{}, \uct{}. 

\item{\bu} $L_\mu\,$, ~$\mu\in\bbc$, ~$c=\l\mu$, 
~$q$~ arbitrary, ~$\dim L_{\mu} ~=~ 1$, cf. \nct{}. 

\item{\bu} $L_{\L,N}\,$, ~$N\in\bbn+1$, ~$q=e^{\pi i/N}$,
~$\mu,c\in\bbc$~ arbitrary not satisfying \vgp{a,b}, $c\neq 0$
for $N$ even, ~$\dim L_{\L,N} ~=~ N$, cf. \vct{}. 

\item{\bu} $L_{n,c,N}\,$, ~$n,N\in\bbn$, $n<N$, 
~$q=e^{\pi i/N}$, ~$c\in\bbc$, $c\neq 0$, 
~$\dim L_{n,c,N} ~=~ n$, cf. \fct{}. 

\item{\bu} $ L_{\mu,\tN}\,$, ~$N=2\tN\in 2\bbn$, 
~$q=e^{\pi i/N}$, ~$\mu\in\bbc$, $\mu\neq 0$, ~$c=0$, 
~$\dim L_{\mu,\tN} ~=~ \tN$, cf. \vctt{}. 

\newsubsec{} 
For the algebra ~$\cb_r$~ with $q\in\bbr$, $q\neq 0$~: 

\item{\bu} $L_{n,c}\,$, ~$n\in\bbn$, ~$c\in\bbr$, $c\neq 0$,
~$\dim L_{n,c} ~=~ n$, cf. \fct{}, \uct{}; unitary for ~$q,c>0$. 

\newsubsec{} 
For the algebra ~$\cc$~: 

\item{\bu} $\tL_{n,\eps}\,$, ~$n\in\bbn$, ~$\eps=\pm1$, 
~$q$~ generic, ~$\dim L_{n,c} ~=~ n$, cf. \wct{}, \uctt{}. 

\item{\bu} $\tL_{\L,N}\,$, ~$N\in\bbn+1$, ~$q=e^{\pi i/N}$, 
~$\mu,c\in\bbc$~ related by \crel\ and not satisfying \vgp{a,b},
$c\neq 0$ for $N$ even, ~$\dim L_{\L,N} ~=~ N$, cf. \vct{}. 

\item{\bu} $\tL_{n,\eps,N}\,$, ~$n,N\in\bbn$, $n<N$, 
~$q=e^{\pi i/N}$, ~$\eps=\pm1$, ~$\dim \tL_{n,\eps,N} ~=~ n$, cf.
\wct{}. 

\item{\bu} $\tL_{\eps,\tN}\,$, ~$N = 2\tN\in 2\bbn$, 
~$q=e^{\pi i/N}$, ~$\dim \tL_{\eps,\tN} ~=~ \tN$, cf. \vctt{}. 

\newsubsec{} Of the above irreps only ~$L_{n,c}\,$~ and ~
$\tL_{n,\eps}$~ have classical ~$sl(2),\ su(2) $~ counterparts.
For fixed $n$ for both cases this is the $n$-dimensional HWM of
$sl(2)$ or $su(2)$ with the conjugation $\om$. The latter HWM is
obtained from ~$L_{n,c}\,$, $\tL_{n,\eps}\,$, resp., 
for ~$q,c\to 1$, $q,\eps\to 1$, resp.

\newsubsec{} 
Of the above irreps all but $L_M\,$, $L_{M,\tN}\,$,
$\tL_{\eps,\tN}\,$ have analogs in the representation theory
\Doa\ of the quantised enveloping algebra $U_q(sl(2))$. However,
the matrix elements there are given by expressions different from
ours. 

\np 

\nt 
{\bf Acknowledgments.}~~ This work was supported by a Royal
Society Grant. 

\vskip 10mm 

\parskip=0pt 
\listrefs 
\np \end